\documentclass[smallextended,numbook,runningheads]{svjour3}
\usepackage[T1]{fontenc}
\usepackage{amsmath,amscd,amssymb}
\usepackage{graphicx}
\usepackage[ruled,lined,linesnumbered]{algorithm2e}
\usepackage{psfrag}

\usepackage{mathptmx}
\journalname{BIT}
\newtheorem{condition}[theorem]{Condition}

\author{ Tomas Johnson \and Warwick Tucker}
\institute{T. Johnson \and W. Tucker \at Department of Mathematics, Uppsala University, Box 480, 751 06 Uppsala, Sweden, \email{tomas.johnson@math.uu.se, warwick.tucker@math.uu.se} } 

\begin{document}
\title{On a computer-aided approach \\ to the computation of Abelian integrals}
\maketitle

\begin{abstract}
An accurate method to compute enclosures of Abelian integrals is developed. This allows for an accurate description of the phase portraits of planar polynomial systems that are perturbations of Hamiltonian systems. As an example, it is applied to the study of bifurcations of limit cycles arising from a cubic perturbation of an elliptic Hamiltonian of degree four. 
\end{abstract}


\section{Introduction}
Nonlinear ordinary differential equations are one of the most common models used in any application of mathematical modelling. 
In this paper we study families of such equations
\begin{equation}\label{ND}
 \left\{ \begin{array}{ccccc}
 \dot{x} & = & -H_y + \epsilon f(x,y) \\
 \dot{y} & = & H_x + \epsilon g(x,y),\end{array} \right.
\end{equation}
depending on a small parameter $\epsilon$.

A fundamental question about such systems is to determine the number and location of limit cycles bifurcating from it as $\epsilon \rightarrow 0$.

In general, the question about the maximal number of limit cycles, and their location, of a polynomial planar vector field is the second part of Hilbert's 16th problem, which is unsolved even for polynomials of degree 2. For an overview of the progress that has been made to solve this problem we refer to \cite{I02}. Results for the degree 2 case, and a general introduction to the bifurcation theory of planar polynomial vector fields can be found in \cite{R98}. What is known, is that any given polynomial vector field can have only a finite number of limit cycles; this is proved in \cite{E92,I91}.

A restricted version of Hilbert's 16th problem, known as the \textit{weak}, or sometimes the \textit{tangential}, or the \textit{infinitesimal}, Hilbert's 16th problem, asks for the number of limit cycles that can bifurcate from a perturbation of a Hamiltonian system, see e.g. \cite{CL07}. The weak Hilbert's 16th problem has been solved for the degree 2 case, see \cite{CLLZ06}.

Special cases of Hamiltonian systems are those coming from a one dimensional system, $H(x,\dot x)=\frac{\dot x^2}{2}+h(x)$, which we study in the example given in Section \ref{eHam}. If one, in addition, assumes that $f=0$, and $g(x,\dot{x})=g(x)\dot{x}$, (\ref{ND}) is known as a Lienard equation. Such equations have been thoroughly studied, and the case where $dH$, and $g$ have degree 3 has been solved, see \cite{DL01a,DL01b,DL03a,DL03b}. We study general $g$ of degree 3; the set-up of the problem is given in Section \ref{eHam}.

In this paper we present a rigorous, computer-aided approach to find limit cycles of planar polynomial vector fields. A different computer-aided approach was introduced by Malo in his PhD-thesis \cite{M94}, (also described in \cite{G95,GM96}) which is based on the concept of a rotated vector field, as introduced in \cite{D53}. Our approach is completely different: we develop a method to rigorously compute what is known as an Abelian integral. A  brief introduction to Abelian integrals is included in Section \ref{AIsection}. The concept of a computer-aided proof in analysis is based on techniques to rigorously enclose the result of a numerical computation. A basis for such a procedure is interval analysis, introduced by Moore in \cite{Mo66}. By calculating with sets rather than floating points, it is possible to obtain guaranteed results on a computer, enabling automated proofs for continuous problems.

We emphasize that the methods developed in this paper are neither restricted to any specific degree of the polynomial functions $f$, and $g$, nor to the structure of the polynomial Hamiltonian $H$. It can be used to compute Abelian integrals of any polynomial perturbation from any family of compact level curves, ovals, of a polynomial Hamiltonian. The method can be used as a computational tool to accurately describe the phase portraits of a family of planar systems. In the example given in this paper, however, we restrict to the case when $f=0$, and $dH$ and $g$ have degree 3. The method also works for integrable, but non-Hamiltonian, planar polynomial systems. For such systems all formulae need to adjusted to include the integrating factor. 

\section{Abelian integrals}\label{AIsection}
A classical method to prove the existence of limit cycles bifurcating from a family of ovals of a Hamiltonian, $\Gamma_h\subset H^{-1}(h)$, depending continuously on $h$, is to study Abelian integrals, or, more generally, the Melnikov function, see e.g. \cite{CL07,GH83}. Some caution, however, must be taken regarding the correspondence between limit cycles and Abelian integrals, see e.g. \cite{DR06}. Given a Hamiltonian system and a perturbation,
\begin{equation}\label{pHDE}
 \left\{ \begin{array}{ccccc}
 \dot{x} & = & -H_y + \epsilon f(x,y) \\
 \dot{y} & = & H_x + \epsilon g(x,y),\end{array} \right.
\end{equation}
the Abelian integral is defined as
\begin{equation}\label{AI}
I(h)=\int_{\Gamma_h} f(x,y)\,dy-g(x,y)\,dx.
\end{equation}
In this paper all systems and perturbations are polynomial. The most important property of Abelian integrals is described by the Poincar\'e-Pontryagin theorem. 
\begin{theorem}[Poincar\'e-Pontryagin]
Let $P$ be the return map defined on some section transversal to the ovals of $H$, parametrised by the values $h$ of $H$, where $h$ is taken from some bounded interval $(a,b)$. Let $d(h)=P(h)-h$ be the displacement function. Then, $d(h)=\epsilon(I(h)+\epsilon\phi(h,\epsilon)), \quad \quad\rm{as} \quad \epsilon \rightarrow 0,$
where $\phi(h,\epsilon)$ is analytic and uniformly bounded on a compact neighbourhood of $\epsilon=0,\, \,h\in(a,b)$.
\end{theorem}
\begin{proof}
see e.g. \cite{CL07}.
\end{proof}
%

\section{Computer-aided computation of Abelian integrals}
\subsection{Computer-aided proofs}
To prove mathematical statements on a computer, we need an arithmetic which gives guaranteed results. Many computer-aided proofs, including the results in this paper, are based on interval analysis, e.g. \cite{GMT03,H05,T02}. Interval analysis yields rigorous results for continuous problems, taking both discretisation and rounding errors into account. For a thorough introduction to interval analysis we refer to \cite{AH83,Mo66,Mo79,Ne90,Pe98}.
\subsection{Outline of the approach}
The main idea of this paper is to develop a very accurate, validated method to enclose the value of a general Abelian integral. Such a method enables us to sample values of $I(h)$. If we can find two ovals $\Gamma_{h_1}$, and $\Gamma_{h_2}$, such that 
\begin{equation}
I(h_1)I(h_2)<0, 
\end{equation}
then there exists $h^*\in (h_1,h_2),$ such that $I(h^*)=0$.

Since $P_\epsilon$, the return map of the perturbed vector field, is analytic and non-constant, it has isolated fixed points. Thus, we have proved the existence of (at least) one limit cycle bifurcating from $\Gamma_{h^*}$.

\subsection{Computing the integrals}\label{compInt}
To compute the Abelian integral (\ref{AI}) of the form (\ref{1F}), we apply Stokes theorem to get 
\begin{equation}
I(h)=\int_{D_h} \,dw,
\end{equation}
where $D_h$ denotes the interior of an oval $\Gamma_h$. The reason why we prefer to calculate surface integrals, rather than contour integrals, is that we cannot represent the ovals of $H$ exactly. We can only find a cover of the ovals, and the area of this cover yields the uncertainty of our calculations, automatically handled by the interval arithmetic. If we had chosen to compute contour integrals, all of our computations would have been subjected to those errors, since we would always integrate over an unknown location. When calculating surface integrals, however, the effect of the uncertainty of the location of the ovals only contributes on a very small portion of the total area of $D_h$. Note that, inside $D_h$ it is possible to integrate $d\omega$ exactly, that is, there are no truncation errors.

The actual computation of the integrals is performed in four steps; first we find a trapping region for the interesting family of ovals, second we adaptively split this region into three parts, one that covers the oval, one representing the inside and one representing the outside, third we change the coordinates on the boxes covering the oval in order to minimise the area of the cover, fourth we integrate $d\omega$ on the boxes representing the inside and the cover of the oval.

The first step is simple, since we primarily study ovals that are situated inside of a homo- or heteroclinic orbit, exterior ones are only studied after choosing the perturbation. A short branch--and--bound algorithm quickly finds a box enclosing the homo- or heteroclinic orbit, and its interior; this box is our initial domain used for the main part of the program.

In the second step -- the adaptive splitting of the domain -- we perform a series of tests to determine whether a box $B$ intersects the oval, is inside it, or outside it. We start by evaluating the Hamiltonian on $B$ using monotonicity and central forms; since the Hamiltonians we study are sufficiently simple, we implement the derivatives symbolically. By monotonicity we mean that if the partial derivatives are non-zero, then an enclosure of $H$ on an entire box is given by the hull of the enclosures of the values of $H$ on the endpoints. A central form for $H$ on a box $B$, with $(x,y)\in B$ is given by:
$$
H(x,y)+H_x(B_1,B_2)(B_1-x)+H_y(B_1,B_2)(B_2-y).
$$
For a given box, $H$ is evaluate three times; naively, with monotonicity, and using a central forms. Finally, all three enclosures are intersected. Three cases occur: if $H<h$, then $B$ is inside the oval, and we label $B$ as such. If $H>h$ then $B$ is outside the oval and we ignore it. Finally, if $h\in H$, then we try to perform the change of variables as described below. If the change of variables procedure fails, and the size of $B$ is greater than some stopping tolerance, {\tt minsize}, then we split the box $B$ into four parts and re--examine them separately. If the size of $B$ is smaller than {\tt minsize}, then it is labelled {\tt fail}. If the change of variables procedure works, then we label $B$ as {\tt on}.

\begin{figure}[h]
\begin{center}
\includegraphics[width=0.45\textwidth]{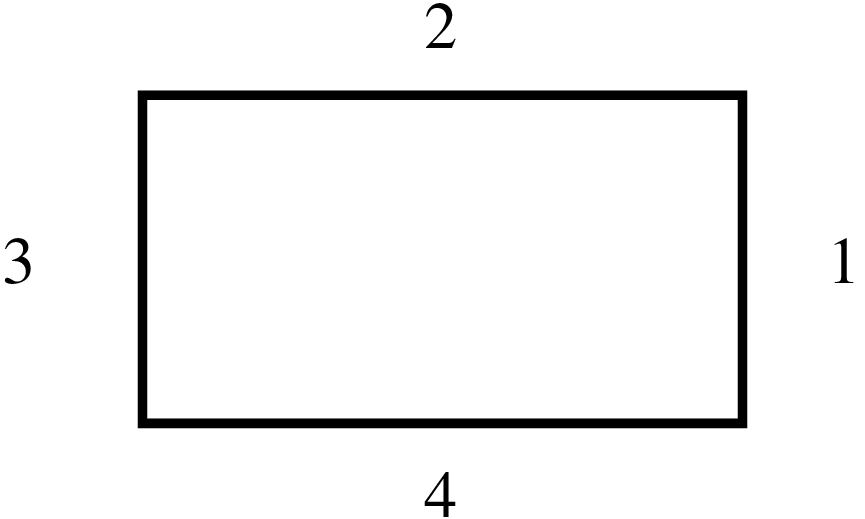}
\caption{The labelling of boxes intersecting an oval.}\label{boxnumb}
\end{center}
 \end{figure}

\begin{figure}[h]
\begin{center}
\includegraphics[width=0.20\textwidth]{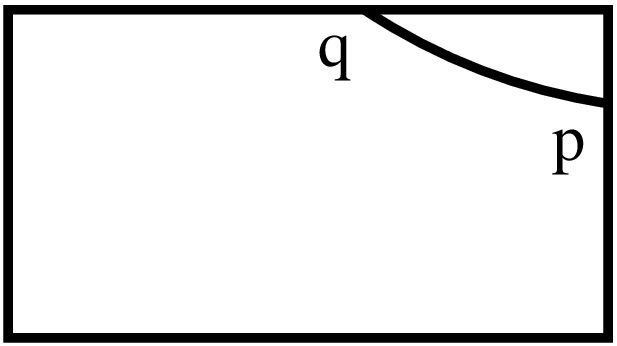}
\includegraphics[width=0.20\textwidth]{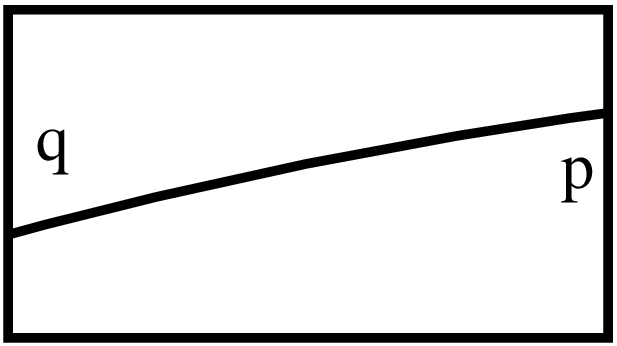}
\includegraphics[width=0.20\textwidth]{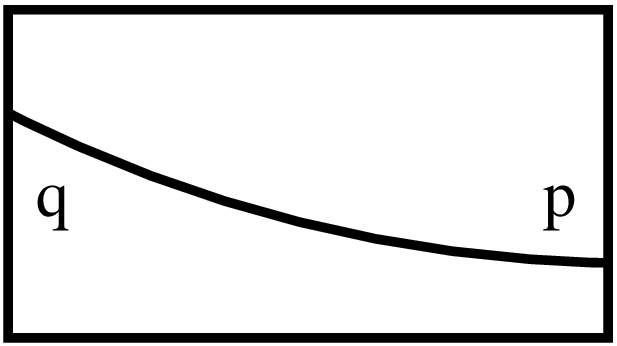}
\includegraphics[width=0.20\textwidth]{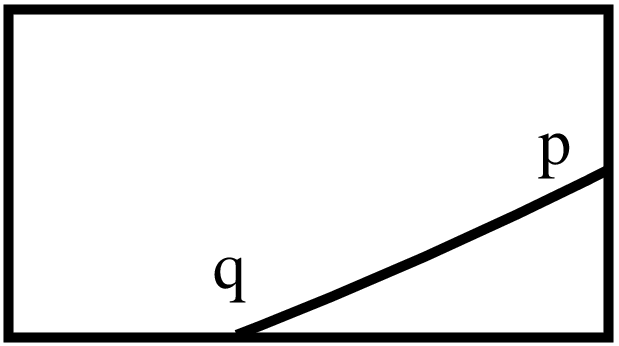}
\includegraphics[width=0.20\textwidth]{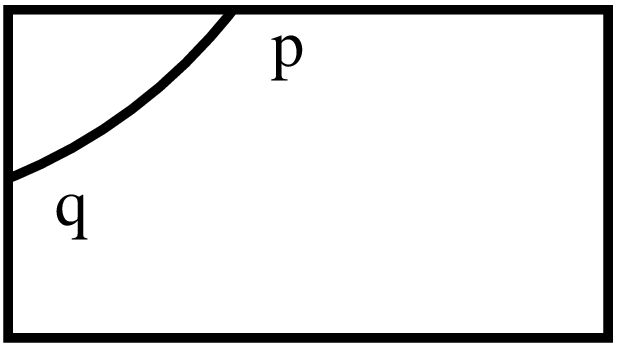}
\includegraphics[width=0.20\textwidth]{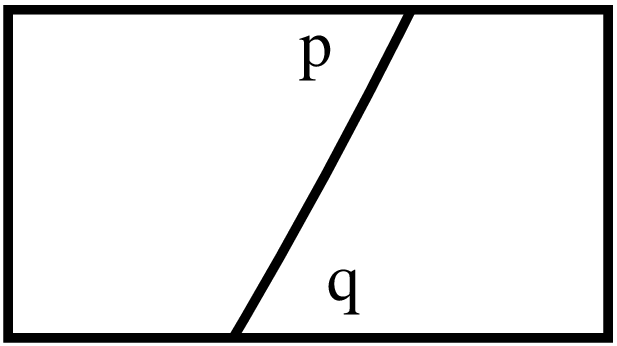}
\includegraphics[width=0.20\textwidth]{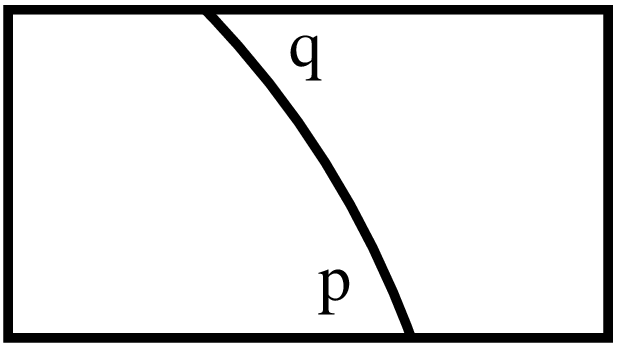}
\includegraphics[width=0.20\textwidth]{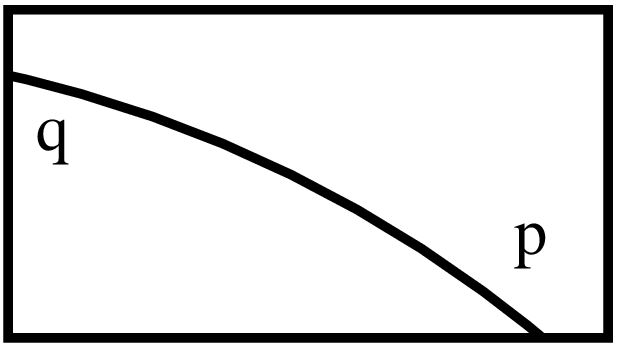}
\caption{The allowed configurations of the intersection of an oval and a box.}\label{config}
\end{center}
\end{figure}

The third, and most complicated, part of our program is the change of variables in the boxes that intersect the oval. Let $b\in B$ be the midpoint of $B$.
Compute $$u=\nabla H(b),$$ and choose $v$ such that $$u\perp v\quad \textrm{and} \quad v_1\geq0.$$ Using the labelling illustrated in Figure \ref{boxnumb}, let {\tt right} and {\tt left} be the sides intersected by the line $b+tv$, $t \in \mathbb{R}$. Denote the intersection points of the straight line $b+tv$ with the boundary of the box by $p$, and $q$, respectively. Note that {\tt right} and {\tt left} are different, since they are the intersections of the boundary of the box with the straight line $b+tv$. We also remark that the possible values of {\tt right} and {\tt left} are {\tt right}$\in\{1,2,4\}$ and {\tt left}$\in \{2,3,4\}$. The allowed configurations of an intersection of the oval with a box are illustrated in Figure \ref{config}. The restriction of $H$ to the sides {\tt right} and {\tt left}, respectively, are one-dimensional functions, and the location of the intersections can be approximated, and their uniqueness proved, using the interval Newton method \cite{Mo66} initialised from the points $p$, and $q$, respectively. If the geometry is not as in Figure \ref{config}, e.g. if the real intersections are on the same side, then uniqueness will fail, and the box $B$ is split, if it is larger than {\tt minsize}, and re-examined.

\begin{figure}[h]
\psfrag{q up}{$q_{up}$}
\psfrag{q down}{$q_{down}$}
\psfrag{p up}{$p_{up}$}
\psfrag{p down}{$p_{down}$}
\psfrag{p}{$p$}
\psfrag{q}{$q$}
\begin{center}
\includegraphics[width=0.65\textwidth]{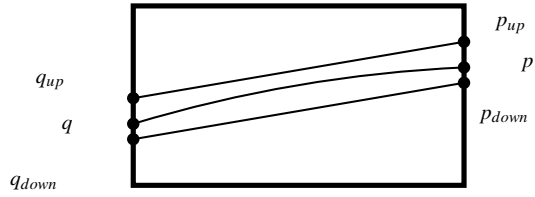}
\caption{Constructing a small, local enclosure of the oval. The geometry corresponds to the second case in Figure \ref{config}.}\label{boxpert}
\end{center}
\end{figure}

Let, $$\tt accuracy=minsize/10.$$ Define the points $p_{up}, \,p_{down}$ on the {\tt right}-side and the points $q_{up}, \,q_{down}$ on the {\tt left}-side at the distance {\tt accuracy} from $p$ and $q$, respectively, as illustrated in Figure \ref{boxpert} for the second case of Figure \ref{config}.

If the following conditions hold, then the oval is inside the tube illustrated in Figure \ref{boxpert}, and we can change coordinates to get a small box, which is guaranteed to contain the segment of the oval passing through $B$. This small box represents the error caused by the unknown location of the oval.

\begin{condition}\label{CV}
$$\begin{array}{ccccc}
& &\begin{rm} sign \end{rm}\left(H(p_{up})-h\right) & = & \begin{rm} sign \end{rm}(H(q_{up})-h) \\& & & = &\begin{rm} -sign \end{rm}(H(p_{down})-h) \\ & & & = & \begin{rm} -sign \end{rm}(H(q_{down})-h)
\end{array}. $$
\end{condition}
Let $l_{up}$, and $l_{down}$ denote the line segments between $p_{up}$ and $q_{up}$, and $p_{down}$ and $q_{down}$, respectively. Denote by $H'$ differentiation with respect to the parametrisation of the line $l_{up}$, and $l_{down}$, respectively.
 
\begin{condition}\label{noCross}
$$ 0 \notin (H(l_{up})-h),$$
and
$$ 0 \notin (H(l_{down})-h).$$
\end{condition}
Let {\tt otherside1}, and {\tt otherside2} be the two other sides of the box $B$, that is, 
$${\tt otherside1}\cup{\tt otherside2}\cup{\tt right}\cup{\tt left}=\{1,2,3,4\}.$$
\begin{condition}\label{noTurn}
$$\Gamma_h \cap {\tt otherside1}=\emptyset  \quad \begin{rm} and \end{rm} \quad \Gamma_h \cap {\tt otherside2}=\emptyset,$$
\end{condition}
Condition \ref{noTurn} is proved using the interval Newton method for the function $H-h$, restricted to {\tt otherside1} and {\tt otherside2}, respectively.

We enclose the segment of the oval inside of the box between two straight lines: Condition \ref{CV} guarantees that the points $p_{up},\,p_{down},\,q_{up},$ and $q_{down}$ are on different sides of the oval as in Figure \ref{boxpert}, Condition \ref{noCross} guarantees that the lines $l_{up}$ and $l_{down}$ do not intersect the oval, and Condition \ref{noTurn} guarantees that the oval does not cross the other sides of the box. Recall that the uniqueness of $p$ and $q$ is proved as they are approximated. Hence, we have proved that the segment of the oval crossing the box has exactly two intersections with the boundary of the box, and that it is confined to the region between $l_{up}$ and $l_{down}$.

If (\ref{CV}), (\ref{noCross}), and (\ref{noTurn}), hold, then we set {\tt accuracy=accuracy/2}, re-calculate $p_{up}$, $p_{down}$, $q_{up},$ and $q_{down}$, and try to verify (\ref{CV}), (\ref{noCross}), and (\ref{noTurn}). This procedure is iterated until (\ref{CV}), or (\ref{noCross}) do not hold. Finally, we label $B$ as {\tt on}.

The fourth and final part of our integration algorithm, is the actual integration. The integration is done separately for the boxes that are labelled, {\tt inside}, {\tt fail}, and {\tt on}.

If $B$ is {\tt inside} we compute 
\begin{eqnarray}\label{IntBox}
\int_B x^iy^j\,dx\wedge dy&=&\left(\frac{\sup(B_1)^{i+1}}{i+1}-\frac{\inf(B_1)^{i+1}}{i+1}\right)\\ \nonumber &\times &\left(\frac{\sup(B_2)^{j+1}}{j+1}-\frac{\inf(B_2)^{j+1}}{j+1}\right).
\end{eqnarray}
If $B$ is labelled {\tt fail}, we know that $B$ might intersect the oval, that is, we have neither been able to prove intersection, nor non-intersection. Therefore, we must include any possible result; the integral over $B$ is calculated as the interval hull of $0$ and the largest, and smallest, respectively, result of (\ref{IntBox}) calculated on a subbox $\tilde B \subset B$. Note that for boxes which intersect the $x$-axis or $y$-axis, this implies that the minus sign between the terms in (\ref{IntBox}) is replaced by a plus sign. E.g., if a box, $B$, is such that $0\in B_1$ and $i=2k-1$, then (\ref{IntBox}) is replaced by:

\begin{eqnarray}\label{IntBox0}
\int_B x^{2k-1}y^j\,dx\wedge dy&=&\pm\left(\frac{\sup(B_1)^{2k}}{2k}+\frac{\inf(B_1)^{2k}}{2k}\right)\\ \nonumber &\times &\left(\frac{\sup(B_2)^{j+1}}{j+1}-\frac{\inf(B_2)^{j+1}}{j+1}\right).
\end{eqnarray}

Boxes labelled {\tt fail} cause large over-estimations. Fortunately such boxes are rare, typically less than $5\%$ of the {\tt on}-boxes, see Section \ref{Res}. If {\tt minsize} is taken sufficiently small, the effect of the {\tt fail}-boxes is negligible.

\begin{figure}[h]
\psfrag{B1}{$B_l$}
\psfrag{B2}{$B_u$}
\psfrag{T1}{$T_l$}
\psfrag{T2}{$T_u$}
\psfrag{P}{$P$}
\begin{center}
\includegraphics[width=0.45\textwidth]{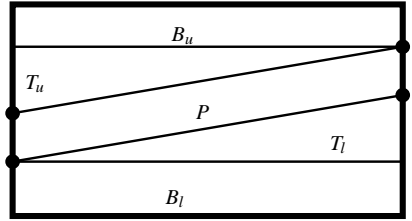}
\caption{The change of variables splitting. The geometry corresponds to the second case in Figure \ref{config}.}\label{split}
\end{center}
\end{figure}

The boxes that are labelled {\tt on}, are split into five parts, as illustrated in Figure \ref{split}. By construction, none of the triangles, $T_l, T_u$, or boxes $B_l, B_u$ in the splitting of $B$ intersect the oval, thus it suffices to evaluate $H$ in one point of each, and hence they can all be labelled as {\tt inside} or {\tt outside}. The boxes $B_l, B_u$ are then treated as above, that is, if they are labelled {\tt inside} they are integrated according to (\ref{IntBox}), and if they are labelled {\tt outside} they are neglected. A triangle labelled {\tt outside} is also neglected, the integrals on triangles labelled {\tt inside} are enclosed by the formula
\begin{equation}\label{IntTri}
\int_T x^iy^j\,dx\wedge dy \in \Box T_1^i\Box T_2^j|T|, 
\end{equation}
where $\Box T$ is the box hull of $T$, and $|T|$ is the area of $T$. This gives a reasonably narrow enclosure of the integral, since the width of $\Box T$ is typically small.
The parallelepiped, $P$, which covers the segment of the oval, remains to be studied. When we integrate over $P$, the same problem as in the {\tt fail} case occurs; we do not know how much of the parallelepiped to include. Therefore, we have to take the hull of all possible outcomes. Hence, the integrals are computed as
\begin{equation}\label{IntPar}
 \int_P x^iy^j\,dx\wedge dy \in \textrm{Hull} \left(0, \Box P_1^i \Box P_2^j|P|\right),
\end{equation}
where $\Box P$ is the box hull of $P$ and $|P|$ is the area of $P$.

The value of the Abelian integral is enclosed by summing over all the computed integrals that are labelled as either {\tt inside}, {\tt fail}, or {\tt on}.
\begin{equation}\label{cAI}
\begin{array}{ccl}
I(h) & \in & \sum_{B\in {\tt inside}} \,(\ref{IntBox}) + \sum_{T\in {\tt inside}} \,(\ref{IntTri}) \\
 & + & \sum_{B\in {\tt fail}} \,\textrm{Hull}(0, (\ref{IntBox})) \\
 & + & \sum_{P\in {\tt on}} \,(\ref{IntPar})
\end{array}
\end{equation}

Thus, we have proved the following:
\begin{theorem}
If Condition \ref{CV}, Condition \ref{noCross}, and Condition \ref{noTurn} hold, then the value of the Abelian integral 
$$ I_{ij}(h) = \int_{D_h} x^iy^j\,dx\wedge dy,$$
is enclosed by Equation \ref{cAI}.
\end{theorem}

The algorithm is given as Algorithm \ref{mainAlgorithm}.
\begin{small}
\begin{algorithm}[h]\label{mainAlgorithm}
 \KwData{$H$, $B$, $h$, $i$, $j$, $minsize$, $accuracy$}
\KwResult{$I_{ij}(h)$}
$I_{ij}(h)=0$\;
workStack+=B\;
\While{notEmtpy(workStack)}{
B=Pop(workStack)\;
\eIf{$H(B)<h$}
{
	$I_{ij}(h)+=\left(\frac{\sup(B_1)^{i+1}}{i+1}-\frac{\inf(B_1)^{i+1}}{i+1}\right)\times \left(\frac{\sup(B_2)^{j+1}}{j+1}-\frac{\inf(B_2)^{j+1}}{j+1}\right)$\;
}
{
	\If{$h\in H(B)$}
	{
		\eIf{Condition \ref{CV} \& Condition \ref{noCross} \& Condition \ref{noTurn}}
		{
			$I_{ij}(h)+=\left(\frac{\sup(B^{l\vee u}_1)^{i+1}}{i+1}-\frac{\inf(B^{l\vee u}_1)^{i+1}}{i+1}\right)\times \left(\frac{\sup(B^{l\vee u}_2)^{j+1}}{j+1}-\frac{\inf(B^{l\vee u}_2)^{j+1}}{j+1}\right) 
			+\Box (T^{l\vee u}_1)^i\Box (T^{l\vee u}_2)^j|T^{l\vee u}|
			+\textrm{Hull} \left(0, \Box P_1^i \Box P_2^j|P|\right)$
		}
		{
			\eIf{diam(B)<minsize}
			{
				\eIf{$0\in B_1$ \& $i$ odd}
				{
					$xInt=\textrm{Hull}\left(-\left(\frac{\sup(B_1)^{i+1}}{i+1}+\frac{\inf(B_1)^{i+1}}{i+1}\right),\left(\frac{\sup(B_1)^{i+1}}{i+1}+\frac{\inf(B_1)^{i+1}}{i+1}\right)\right)$\;
				}
				{
					$xInt=\left(\frac{\sup(B_1)^{i+1}}{i+1}-\frac{\inf(B_1)^{i+1}}{i+1}\right)$\;
				}
				\eIf{$0\in B_2$ \& $j$ odd}
				{
					$yInt=\textrm{Hull}\left(-\left(\frac{\sup(B_2)^{j+1}}{j+1}+\frac{\inf(B_2)^{j+1}}{j+1}\right),\left(\frac{\sup(B_2)^{j+1}}{j+1}+\frac{\inf(B_2)^{j+1}}{j+1}\right)\right)$\;
				}
				{
					$yInt=\left(\frac{\sup(B_2)^{j+1}}{j+1}-\frac{\inf(B_2)^{j+1}}{j+1}\right)$\;
				}
				$I_{ij}(h)+=\textrm{Hull} \left(0,xInt\times yInt\right)$\;
			}
			{
				splitAndStore(B,workStack)\;
			}
		}
	}
}
} 
\caption{Implementation of the Algorithm} 
\end{algorithm}
\end{small}
\section{Computational results}\label{Res}
In this section we apply the methods developed in Section \ref{compInt} to an elliptical Hamiltonian of degree four, described in Section \ref{eHam}. The main idea  is to integrate monomial forms at some points, and then to specify the coefficients of the perturbation $\omega$ such that $I(h)$ is zero at the sampled points. Therefore, let
\begin{equation}
I_{ij}(h)=\int_{D_h} x^iy^j\,dx\wedge dy.
\end{equation}
We sample at some number of $h$-values, uniformly distributed between the saddle loops and the singularity. From these calculations we deduce candidate coefficients.

Given some candidate coefficients of the form $\omega$, we calculate the $I_{ij}(h)$, at intermediate ovals. If the linear combination of the $I_{ij}(h)$ has validated sign changes between the sample points we are done: it has been proved that the corresponding perturbation yields bifurcations with the given number of limit cycles as $\epsilon \rightarrow 0$.

All computations were performed on a Intel Xeon 2.0 Ghz, 64bit processor with 7970Mb of RAM. The program was compiled with \texttt{gcc}, version 3.4.6. The software for interval arithmetic was provided by the \texttt{C-XSC} package, version 2.1.1, see \cite{CXSC,HH95}.

\subsection{Example - bifurcations from a figure eight loop}\label{eHam}
\begin{figure}[h]
\begin{center}
\includegraphics[width=0.6\textwidth]{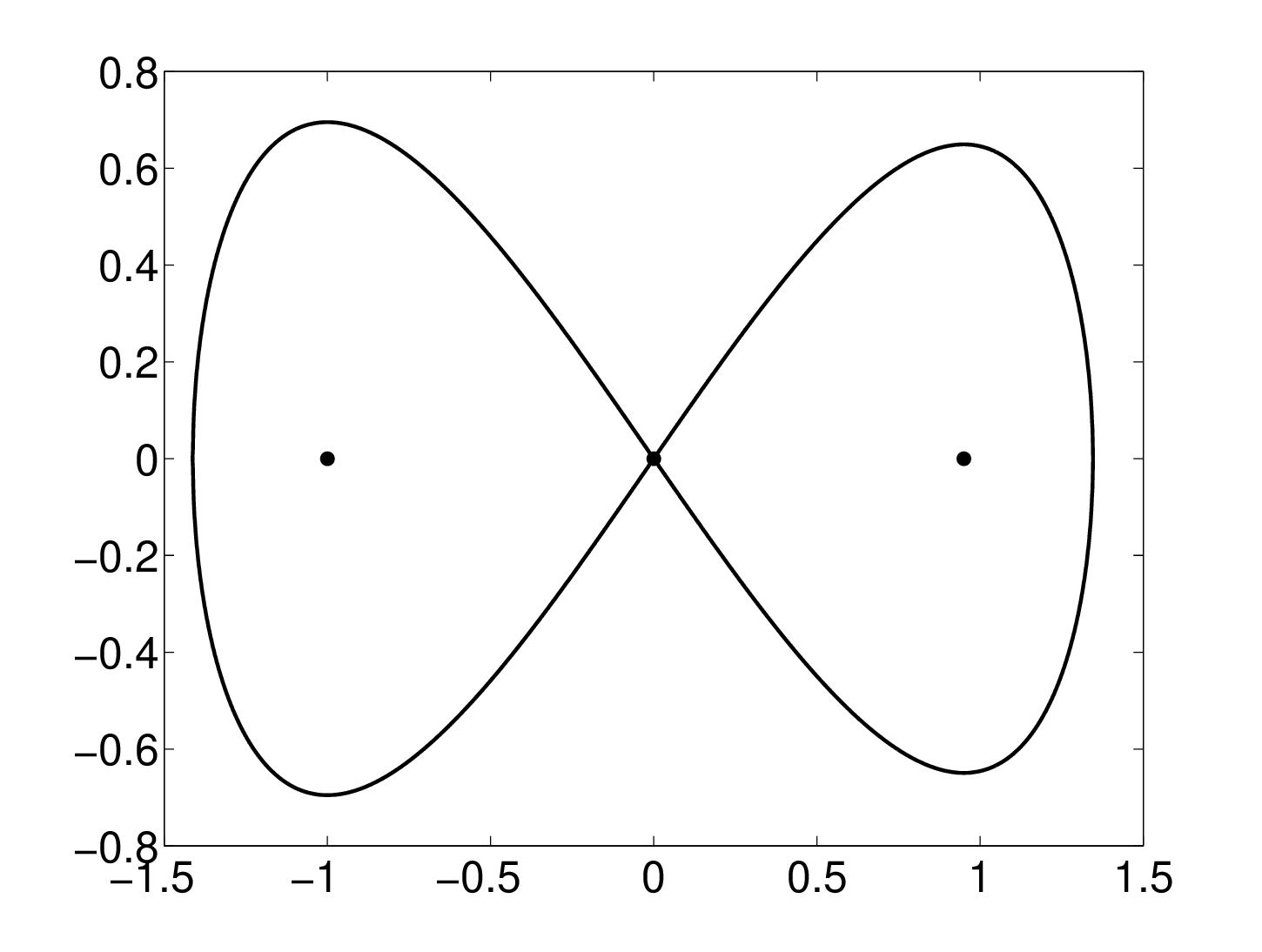}
\caption{The elliptic Hamiltonian of degree 4 with a figure eight loop, studied in section \ref{eHam}.}\label{Hpic}
\end{center}
 \end{figure}

We study the elliptic Hamiltonian of degree 4 with a figure eight loop, given by
\begin{equation}\label{Ham}
 H=\frac{y^2}{2}+\frac{x^4}{4}+\frac{1-\lambda}{3}x^3-\frac{\lambda}{2} x^2,
\end{equation}
where $\lambda\in(0,1)$, see \cite{DL03b}. The corresponding differential system has two centres, at $H=-\frac{1}{12}(2\lambda+1)$, and $H=-\frac{1}{12}\lambda^3(\lambda+2)$, that are  surrounded by a figure eight loop, located at $H=0$, see Figure \ref{Hpic}. As $\lambda$ grows the right loop grows; $\lambda=1$ is a symmetric figure eight loop. We choose to study $\lambda=0.95$; a motivation why we want $\lambda$ large is as follows. We want to construct a nontrivial example with as many limit cycles as possible. In the symmetric case the two branches are identical. Therefore, heuristically, it is a reasonable that for $\lambda$ close to one it should be possible to choose coefficients so that the two branches oscillate together. After some experiments we decide to put $\lambda=0.95$, since it is relatively far away from $1$ to be significantly different, but still sufficiently close to $1$ for the domains of the two branches to have a large overlap. This allows us to locate extra limit cycles, compared to what is possible by simply solving the linear system as described below. 

The Hamiltonian (\ref{Ham}) corresponds to the differential system,
\begin{equation}\label{DE}
 \left\{ \begin{array}{ccccc}
 \dot{x} & = & -H_y & = & -y \\
 \dot{y} & = & H_x & = & x^3+(1-\lambda)x^2-\lambda x.\end{array} \right.
\end{equation}
We are interested in limit cycles bifurcating from the periodic solutions of (\ref{DE}), corresponding to integral curves of (\ref{Ham}). The closed level-curves of (\ref{Ham}) are called \textit{ovals}. In a series of papers \cite{DL01a,DL01b,DL03a,DL03b}, Dumortier and Li study cubic perturbations of elliptic Hamiltonians corresponding to Lienard equations. That is, 
\begin{equation}\label{LE}
 \ddot{x} +\epsilon(\alpha +\beta x+\gamma x^2)\dot{x} + ax^3+bx^2+cx = 0.
\end{equation}

For the elliptic Hamiltonians of degree four with compact ovals, there are five different classes of phase portraits, see e.g \cite{CL07}. They are, \textit{the truncated pendulum, the saddle loop. the global centre, the cuspidal loop,} and \textit{the figure-eight loop}. 

Compared to the Lienard case, we add a fourth term, $\delta \frac{y^3}{3}$, to the perturbation, and explore what kind of bifurcations we can prove to exist. We study the perturbed system, 
\begin{equation}\label{pDE}
 \left\{ \begin{array}{ccc}
 \dot{x} & = & -y  \\
 \dot{y} & = & x^3+(1-\lambda)x^2-\lambda x + \epsilon\left((\alpha +\beta x+\gamma x^2)y+\delta\frac{y^3}{3}\right).\end{array} \right.
\end{equation}
The 1-form associated with this perturbation is 
\begin{equation}\label{1F}
 \omega = -\left((\alpha +\beta x+\gamma x^2)y+\delta\frac{y^3}{3}\right)\,dx.
\end{equation}
For computational efficiency we primarily study its exterior derivative,
\begin{equation}\label{2F}
 d\omega = \left((\alpha +\beta x+\gamma x^2)+\delta y^2\right)\,dx\wedge dy.
\end{equation}

In \cite{P91} Petrov proves that when restricting to one family of ovals, surrounding one of the two centres, the space of Abelian integrals has dimension 4, and that the space has the Chebyshev property, that is, the number of zeros of a function in this space is less than the dimension of the space. He also proves that this bound is sharp. To construct an example with more than three limit cycles surrounding either of the two centres, we can therefore not simply use the Chebyshev property of the space of Abelian integrals. 

Our heuristic argument to guess parameters is the following: we start by integrating at 100 uniformly distributed ovals, in each eye of the loop. We do this with moderate accuracy, which gives a fast and sufficiently precise result. Since we have chosen to study a figure eight loop that is not far from being symmetric, it is reasonable to assume that the two branches behave similarly, which makes it probable that we should be able to determine coefficients so that each branch has two zeros. To determine such zeros, we solve the following linear system:

\begin{footnotesize}
\begin{equation}
\left[
\begin{array}{cccc}
I^l_{00}(-0.0362) & I^l_{10}(-0.0362) & I^l_{20}(-0.0362) & I^l_{02}(-0.0362) \\
I^l_{00}(-0.1208) & I^l_{10}(-0.1208) & I^l_{20}(-0.1208) & I^l_{02}(-0.1208) \\
I^l_{00}(-0.1812) & I^l_{10}(-0.1812) & I^l_{20}(-0.1812) & I^l_{02}(-0.1812) \\
I^r_{00}(-0.1054) & I^r_{10}(-0.1054) & I^r_{20}(-0.1054) & I^r_{02}(-0.1054)
\end{array} \right]
\left[
\begin{array}{c}
\alpha \\ \beta \\ \gamma \\ \delta
\end{array} \right]
=\left[
\begin{array}{c}
1 \\ -1 \\ 1 \\ -1
\end{array} \right]
\end{equation}
\end{footnotesize}
where $I^l_{ij}(h)$, and $I^r_{ij}(h)$, denote the monomial Abelian integrals calculated on the left and right ovals, respectively. 

\begin{figure}[h]
\begin{center}
\includegraphics[width=0.6\textwidth]{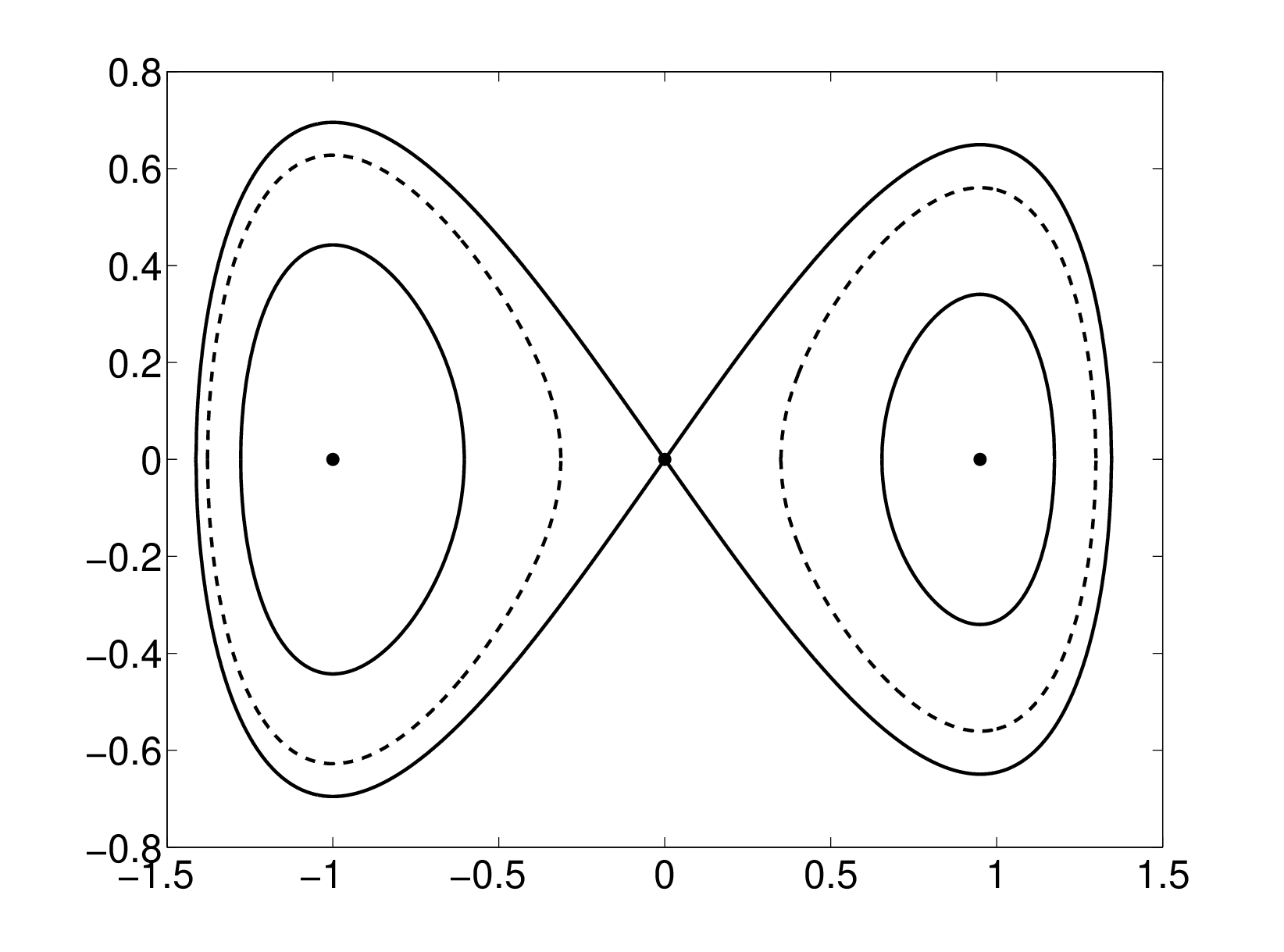}
\caption{The limit cycles from Example \ref{eHam}; unstable limit cycles are dashed. Note that we only prove the existence of the limit cycles, their locations as drawn in the figure are the locations of the ovals they bifurcate from.
}\label{Respic}
\end{center}
\end{figure}

\begin{figure}[h]
\begin{center}
\includegraphics[width=0.75\textwidth]{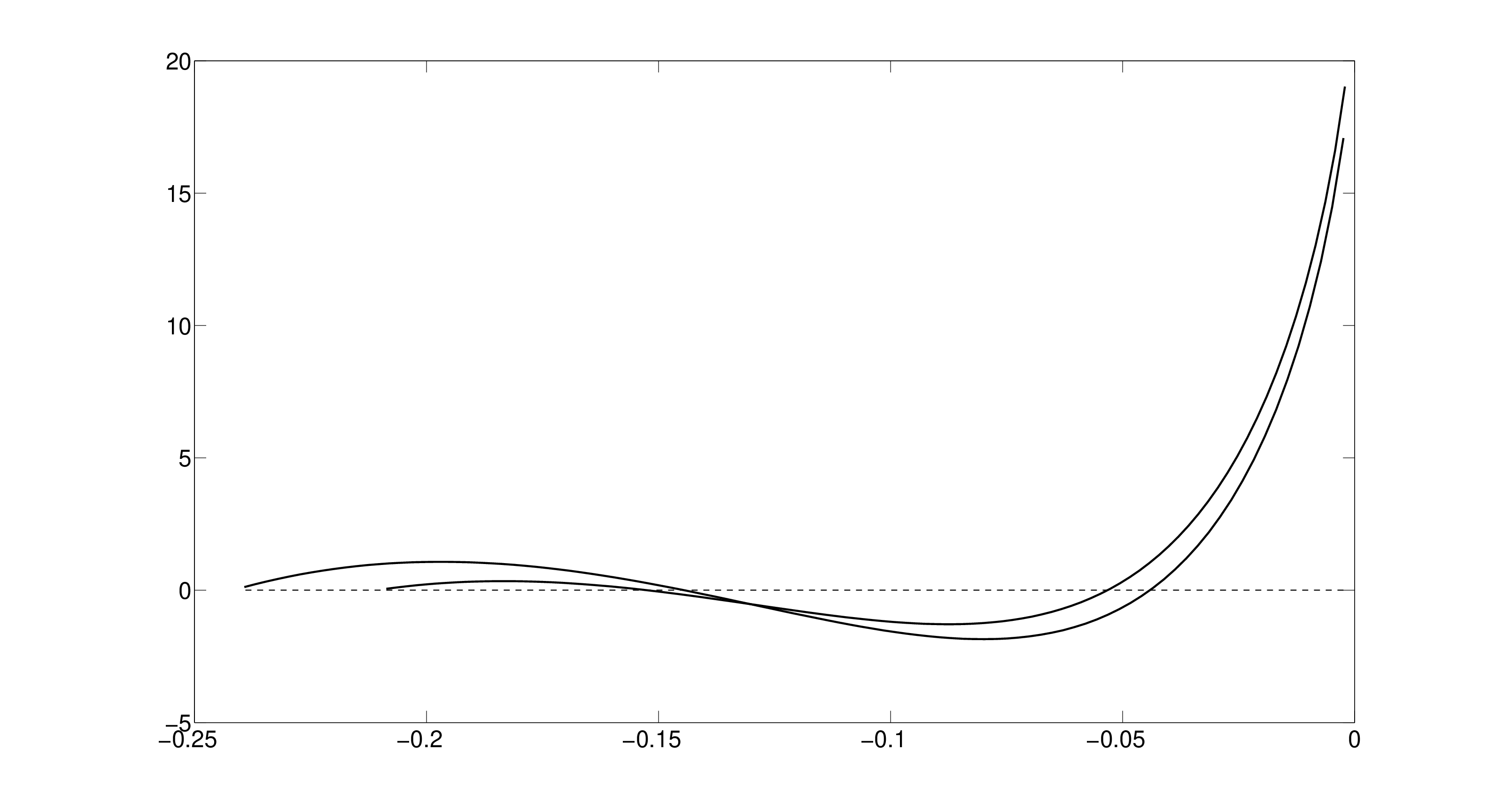}
\caption{The two branches of the Abelian Integral for the figure eight loop.}\label{Ipic8}
\end{center}
\end{figure}

This computation gives the approximate solution $\alpha = 438.4905$, $\beta = -25.2469$, $\gamma = -452.7899$, and $\delta = -741.0341$, which we use as our perturbation. The graph of the resulting function is given in Figure \ref{Ipic8}, which appears to have 4 zeros, illustrated in Figure \ref{Respic}. This, of course, has to be proved.

To prove that the perturbation constructed above has 4 zeros, we procede as in the previous examples, and compute enclosures of the Abelian integral at intermediate ovals. On the left branch we calculate $I(-0.0121)$, $I(-0.0846)$, and $I(-0.1933)$, and on the right branch we compute $I(-0.0105)$, $I(-0.0738)$, and $I(-0.1686)$. The result is given in Tables \ref{8LIntl}, and \ref{8LIntr}.

\begin{table}[h]
\begin{center}
\begin{footnotesize}
\begin{tabular}{l|llll}
$h$ & $I^l_{00}$ & $I^l_{10}$ & $I^l_{20}$ & $I^l_{02}$ \\ \hline
-0.0121 & [1.206,1.207] & [-1.034,-1.033] & [0.9945,0.9951] & [0.1290,0.1293] \\
-0.0846 & [0.7661,0.7665] & [-0.7073,-0.7068] & [0.6902,0.6907] & [0.05829,0.05839] \\
-0.1933 & [0.2219,0.2222] & [-0.2178,-0.2175] & [0.2160,0.2164] & [0.00532,0.00534]
\end{tabular}
\caption{The computed enclosures for the left branch of the figure eight loop.}\label{8LIntl}
\end{footnotesize}
\end{center}

\begin{center}
\begin{footnotesize}
\begin{tabular}{l|llll}
$h$ & $I^r_{00}$ & $I^r_{10}$ & $I^r_{20}$ & $I^r_{02}$ \\ \hline
-0.0105 & [1.077,1.078] & [0.8773,0.8778] & [0.8033,0.8039] & [0.1006,0.1008] \\
-0.0738 & [0.6846,0.6851] & [0.6002,0.6006] & [0.5573,0.5577] & [0.04545,0.04553] \\
-0.1686 & [0.1984,0.1987] & [0.1848,0.1850] & [0.1744,0.1747] & [0.004154,0.004164]
\end{tabular}
\end{footnotesize}
\end{center}
\caption{The computed enclosures for the right branch of the figure eight loop.}\label{8LIntr}
\end{table}

\begin{figure}[h]
\begin{center}
\includegraphics[width=0.9\textwidth]{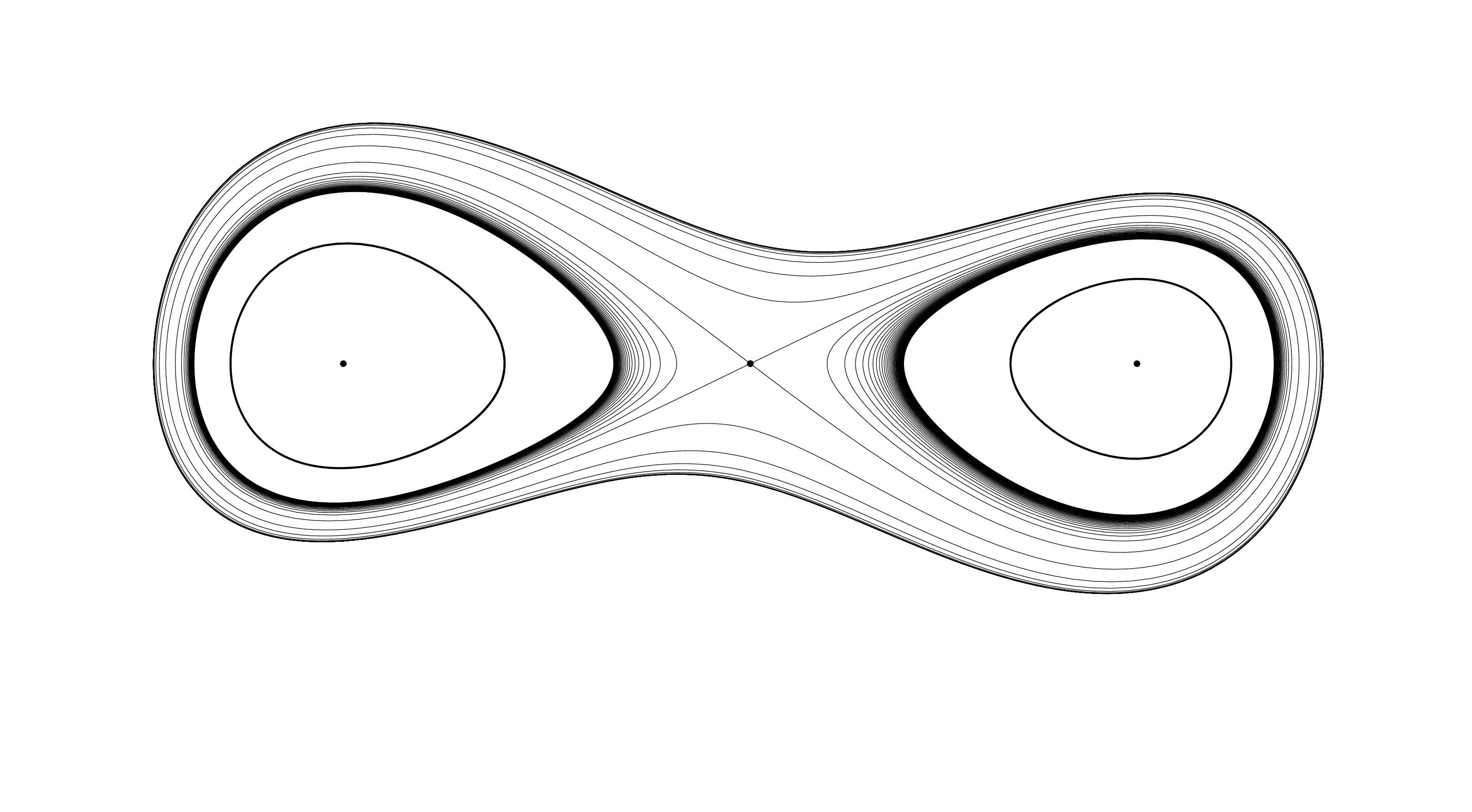}
\caption{The perturbed figure eight loop, here with $\epsilon = 0.001$, illustrating the 5 limit cycles found in Section \ref{eHam}.}\label{Case4}
\end{center}
\end{figure}

Finally, we compute $I^l(h)$, and $I^r(h)$ at the intermediate ovals,
\begin{equation}
\begin{array}{ccl}
I^l(-0.0121) & = & [+8.698,+9.290], \\
I^l(-0.0846) & = & [-2.204,-1.780], \\
I^l(-0.1933) & = & [+0.9121,+1.119],  \\
I^r(-0.0105) & = & [+11.56,+12.10],  \\
I^r(-0.0738) & = & [-1.181,-0.7959], \\
I^r(-0.1686) & = & [+0.2095,+0.3847]
\end{array}.
\end{equation}
Hence, the system with the given perturbation has four limit cycles, one attracting and one repelling inside each loop, see Figures \ref{Respic} and \ref{Case4}.
The run-time of the program was, for the left (right) branch, 82 (78)  seconds, a total of 1182 (1166) boxes were used to cover the 3 ovals, 82 (56) of these belong to the {\tt fail} class. 

To prove that the unstable separatrices of the saddle are attracted to a limit cycle enclosing the figure eight loop, as indicated in Figure \ref{Case4}, we first calculate $I^o(h)$, the outer Abelian integral, for some $h>0$ values with low accuracy to find an indication of a sign change. It appears that a limit cycle bifurcates from an oval close to $H=0.1$. Therefore, we compute $I^o(0.09)$, and $I^o(0.11)$, the result is given in Table \ref{8LInto},

\begin{table}
\begin{center}
\begin{footnotesize}
\begin{tabular}{l|llll}
$h$ & $I^o_{00}$ & $I^o_{10}$ & $I^o_{20}$ & $I^o_{02}$ \\ \hline
0.09 & [3.576,3.587] & [-0.1843,-0.1709] & [2.560,2.575] & [0.5307,0.5376] \\
0.11 & [3.776,3.786] & [-0.1862,-0.1740] & [2.708,2.724] & [0.6044,0.6109] 
\end{tabular}
\end{footnotesize}
\end{center}
\caption{The computed enclosures for the outside of the figure eight loop.}\label{8LInto}
\end{table}

\begin{equation}
\begin{array}{ccl}
I^o(0.09) & = & [+8.715,+24.83],  \\
I^o(0.11) & = & [-25.37,-9.821].

\end{array}
\end{equation}

These calculations verify that the perturbed system has an attracting limit cycle bifurcating from an oval outside the figure eight loop.
The run-time of the program was 39 seconds, a total of 496 boxes were used to cover the 2 ovals, 52 of these belong to the {\tt fail} class. 

\begin{figure}[h]
\begin{center}
\includegraphics[width=0.9\textwidth]{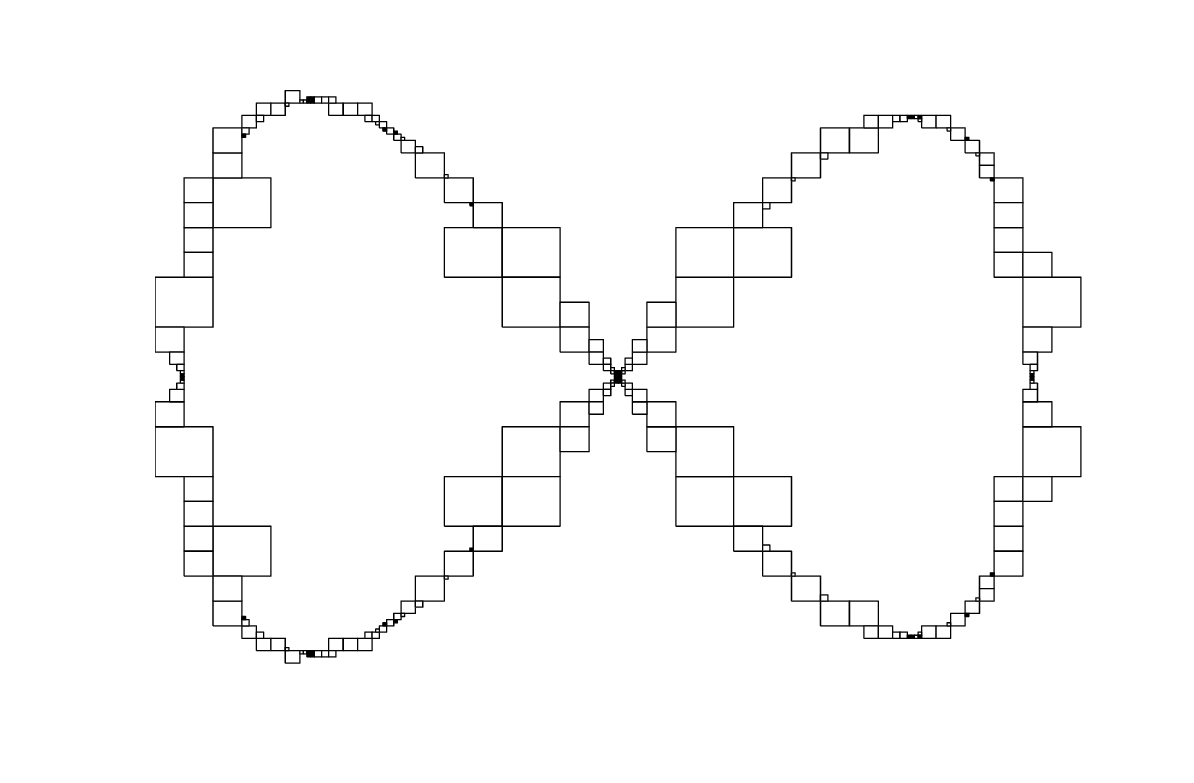}
\caption{A cover of the outer oval with $H=10^{-6}$. The boxes labeled {\tt fail} are plotted as filled black boxes, the boxes that are labeled {\tt on} are plotted as white boxes with a black frame.}\label{cover}
\end{center}
\end{figure}

To illustrate how the algorithm partitions the original trapping-region of an oval into a sufficiently fine cover of it, a cover of an outer oval of the figure eight loop is given in Figure \ref{cover}. Note that the cover is highly non-uniform.

\section{Conclusions}
We have presented a method to rigorously calculate Abelian integrals. The method can be applied to study any polynomial perturbation of a planar polynomial integrable vector field. As an application, we have applied the method to an elliptic Hamiltonian of degree 4.

The method can be used in several ways: either one can use it to verify that a specific perturbation guessed by some other method indeed has a certain number of zeros, or one can use it as in Section \ref{eHam} to sample and plot the monomial Abelian integrals. In the latter case, if a good choice of parameters can be made from the approximate knowledge of the monomial Abelian integrals, then one can re-use the program to verify that guess, as is done in Section \ref{eHam}. This means that one can use the method for experimental, but rigorous, studies of the possible configurations of limit cycles bifurcating from a given planar polynomial Hamiltonian system. We believe that enabling such rigorous studies can be very useful when studying a given system, since the phenomena are typically subtle and hard to detect using floating point computations. Without the verification step it is hard to decide from the computations what is a bifurcation, and what is just numerical noise. 

A major challenge is to device a method which can be used to guess what perturbations to investigate. One such method that appears in the literature is that of a detection function, as used in e.g. \cite{ZXLZ07}. Another problem, which we have ignored in this paper, is that typically when one has a Hamiltonian depending on parameters, the maximal number of limit cycles that can bifurcate from one member of this family, will only appear for some special values of the parameters. It would therefore be desirable to develop conditions indicating how to choose one candidate system from a family. In Section \ref{eHam} we give a completely heuristic argument why we want to have $\lambda$ large. Another method to choose some of the parameters in the Hamiltonian is to restrict the study to systems with maximal number of centres, i.e., of the form

\begin{equation}\label{symHam}
 \left\{ \begin{array}{ccc}
 \dot{x} & = & -y(y^2-b_1)(y^2-b_2)\cdots(y^2-b_k) \\
 \dot{y} & = & x(x^2-a_1)(x^2-a_2)\cdots(x^2-a_k) \end{array} \right.
\end{equation}

\noindent where the $a_i$'s and $b_i$'s are increasing sequences of positive numbers. Different choices of $a_i$ and $b_i$ introduce different symmetries into the system, which can be used to find perturbed systems with a large number of limit cycles.
\begin{small}

\end{small}
\end{document}